\documentclass{amsart}
\usepackage{amsmath,amsthm,amsfonts,amscd,eucal,latexsym}
\usepackage{epsfig}
\usepackage{subfigure}

\numberwithin{equation}{section}

\newcommand{\ch}{{\mathcal H}}
\newcommand{\cai}{{\mathcal I}}
\newcommand{\ct}{{\mathcal T}}

\newcommand{\bn}{{\mathbb N}}
\newcommand{\bq}{{\mathbb Q}}
\newcommand{\br}{{\mathbb R}}

\renewcommand{\d}{\delta}        
\newcommand{\eps}{\varepsilon} 
\renewcommand{\th}{\vartheta}    
\renewcommand{\k}{\kappa}
\renewcommand{\l}{\lambda}       
\newcommand{\La}{\Lambda}
\newcommand{\m}{\mu}
\newcommand{\n}{\nu}
\newcommand{\s}{\sigma}        
\renewcommand{\S}{\Sigma}
\newcommand{\f}{\varphi}
\renewcommand{\O}{\Omega}

\newcommand{\ov}{\overline}
\newcommand{\itm}[1]{\item{$(#1)$}}
\newcommand{\vol}{vol}
\newcommand{\e}[1]{\text{e}^{#1}}

\newcommand{\subc}{\underline{\delta}}
\newcommand{\supc}{\overline{\delta}}
\newcommand{\subd}{\underline{d}}
\newcommand{\supd}{\overline{d}}
\newcommand{\cl}[1]{\text{Cl}(#1)}
\newcommand{\cpt}[1]{\text{Comp}(#1)}
\newcommand{\B}{\ov{B}}
\newcommand{\ceiling}[1]{\lceil #1 \rceil}

\newtheorem{thm}{Theorem}[section]
\newtheorem{cor}[thm]{Corollary}
\newtheorem{prop}[thm]{Proposition}
\newtheorem{lem}[thm]{Lemma}
\theoremstyle{definition}
\newtheorem{dfn}[thm]{Definition}

\newtheorem{cond}[thm]{Condition}
\theoremstyle{remark}
\newtheorem{rmk}[thm]{Remark} 
 
\begin{document}
\title{Tangential dimensions II. Measures}
\author{Daniele Guido, Tommaso Isola}
\address{
Dipartimento di Matematica, Universit\`a di Roma ``Tor
Vergata'', I--00133 Roma, Italy.}
\email{guido@mat.uniroma2.it, isola@mat.uniroma2.it}

\date{\today}
\begin{abstract}
    Notions of (pointwise) tangential dimension are considered, for
    measures of $\br^N$.  Under regularity conditions (volume
    doubling), the upper resp.  lower dimension at a point $x$ of a
    measure $\mu$ can be defined as the supremum, resp.  infimum, of
    local dimensions of the measures tangent to $\mu$ at $x$.  Our
    main purpose is that of introducing a tool which is very sensitive
    to the "multifractal behaviour at a point" of a measure, namely
    which is able to detect the "oscillations" of the dimension at a
    given point, even when the local dimension exists, namely local
    upper and lower dimensions coincide.  These definitions are tested
    on a class of fractals, which we call translation fractals, where
    they can be explicitly calculated for the canonical limit measure. 
    In these cases the tangential dimensions of the limit measure
    coincide with the metric tangential dimensions of the fractal
    defined in \cite{GuIs11}, and they are constant, i.e. do not
    depend on the point.  However, upper and lower dimensions may
    differ.  Moreover, on these fractals, these quantities coincide
    with their noncommutative analogues, defined in previous papers
    \cite{GuIs9,GuIs10}, in the framework of Alain Connes'
    noncommutative geometry.
 \end{abstract}
 \subjclass{28A80; 28A78}
 \keywords{Tangent measure, translation fractals}

\maketitle

 \setcounter{section}{0}

  \section{Introduction.}\label{sec:zeroth}
  
In this paper we continue the analysis concerning notions of
tangential dimesions.  

Our aim is that of finding dimensions describing the non-regularity, or
fractality, of a given measure.  The kind of non-regularity we
study here is related to the fact that a dimension may have an
oscillating behavior at a point.  Indeed dimensions are often defined
as limits, and an oscillating behavior means that the upper and lower
versions of the considered dimension are different.  Our main goal
here is to associate to a measure a local dimension that is able to
maximally detect such an oscillating behavior, namely for which the
upper and lower determinations form a maximal dimensional interval.

Let us recall that we introduced first tangential dimensions in the
framework of Alain Connes' noncommutative geometry \cite{Co}, as
extremal points of the singular traceability interval \cite{GuIs9}.
Their explicit formulas suggested the definition of tangential
dimensions at a point for a metric space, given in  \cite{GuIs11},
where we showed that, under regularity conditions (cf.  Theorem
\ref{newformula} $(i)$ and $(ii)$), upper, resp.  lower, tangential
dimension of a metric space at a given point can be equivalentrly
defined as the supremum, resp.  infimum, of local dimensions of the
tangent sets {\it a la Gromov} at the point.

Here we define tangential dimensions at a point for measures and show
that, under the volume doubling condition, upper, resp.  lower,
tangential dimension of a measure at a given point can be
equivalently defined as the supremum, resp.  infimum, of local
dimensions of the tangent measures at the point.

Finally we compute the tangential dimensions for some classees of
fractals.

Indeed we give a condition (Condition \ref{NewAssum}) on a measure
$\m$ on a metric space $X$ under which the tangential dimensions for
$\m$ coincide with the tangential dimensions of $X$ and are locally
constant.  Furthermore, under the same condition, tangential
dimensions for metric spaces and measures are extrema of local
dimensions of the corresponding tangent objects.

First we consider the class of self-similar fractas with open set 
condition, showing that in this case tangential dimensions do not 
give new information, indeed they coincide with the Hausdorff dimension.

Then we consider the class of translation fractals, show that the
mentioned condition is satisfied by translation fractals with open set
condition for the canonical limit measure and compute the tangential
dimensions for such measure.  Besides their coincidence with metric
tangential dimensions, which follows by the results described above,
direct inspection shows that they also coincide with the
noncommutative tangential dimensions for the spectral triples
associated to translation fractals computed in \cite{GuIs9} and
\cite{GuIs10}.

 \section{Tangential dimensions of measures} 

 In this section we shall define upper and lower tangential dimensions
 of a measure on a metric space $X$ and study some of their
 properties.  The name tangential is motivated by the results in
 subsection \ref{subsub:TangentialFormula}, where we show that for
 Radon measures on $\br^{N}$, under volume doubling condition, the
 upper, resp.  lower, tangential dimension, is simply the supremum,
 resp.  infimum, of the (upper, resp.  lower) local dimensions of the
 tangent measures.
 
\subsection{Basic properties}
 
 Let $(X,d)$ be a metric space, $\m$ a locally finite Borel measure,
 namely $\m$ is finite on bounded sets, and set $B(x,r):= \{y\in X : 
 d(x,y)<r\}$.
 
 Let us recall that the local dimensions of a measure at $x$ are
 defined as
 \begin{align*}
	 \underline{d}_{\m}(x)&=\liminf_{r\to0}
	 \frac{\log\m(B(x,r))}{\log r},\cr
	 \overline{d}_{\m}(x)&=\limsup_{r\to0}
	 \frac{\log\m(B(x,r))}{\log r}.
 \end{align*}
 
 \begin{rmk}
     If $\m$ is zero on a neighborhood of $x$, we set
     $\underline{d}_{\m}(x)= \overline{d}_{\m}(x) =+\infty$.  Indeed,
     let us introduce the following partial order relation on
     measures: $\m<_{x}\n$ if there exists a neighbourhood $\O$ of $x$
     such that for any positive Borel function $\f$ supported in $\O$ we
     have $\langle\m,\f\rangle\leq\langle\n,\f\rangle$.
 
     By definition the maps $\m\mapsto \underline{d}_{\m}(x)$,
     $\m\mapsto \overline{d}_{\m}(x)$ are decreasing, namely reverse
     the ordering.  In particular, if $x$ is not in the support of
     $\m$, namely $\mu$ is zero on a neighbourhood of $x$, the local
     dimensions of $\m$ should be set to $+\infty$.
 \end{rmk}

 Now we introduce tangential dimensions for $\m$.
 
 \begin{dfn}\label{meastgdims}
     The {\it lower and upper tangential dimensions} of $\m$ are defined as
     \begin{align*}
	 \underline{\d}_{\m}(x) & := \liminf_{\l \to 0} \liminf_{r \to 0}
	 \frac{ \log \left( \frac{\m(B(x,r))} {\m(B(x,\l
	 r))} \right) } {\log 1/\l} \in [0,\infty], \\
	 \overline{\d}_{\m}(x) & := \limsup_{\l \to 0} \limsup_{r \to 0}
	 \frac{ \log \left( \frac{\m(B(x,r))}{\m(B(x,\l
	 r))} \right) } {\log 1/\l} \in [0,\infty].
     \end{align*}
 \end{dfn}

 In the following we shall set $f(t)=f_{x,\m}(t) =
 -\log(\m(B(x,e^{-t})))$, and $g(t,h)=f(t+h)-f(t)$.  With this
 notation, the definitions above become
 \begin{align}
     \underline{d}_{\m}(x)&=\liminf_{t\to+\infty}
     \frac{f(t)}{t}
      =  \lim_{t \to \infty}\liminf_{h\to\infty}
     \frac{g(t,h)} {h},\label{eq1}\\
     \overline{d}_{\m}(x)&=\limsup_{t\to+\infty}
     \frac{f(t)}{t}=  \lim_{t \to \infty}\limsup_{h\to\infty}
     \frac{g(t,h)} {h},\label{eq2}\\
     \underline{\d}_{\m}(x) & = \liminf_{h\to\infty} \liminf_{t \to \infty}
     \frac{g(t,h)} {h},  \label{eq3}\\
     \overline{\d}_{\m}(x) & = \limsup_{h\to\infty} \limsup_{t \to \infty}
     \frac{g(t,h)} {h}. \label{eq4}
 \end{align}

 \begin{thm}\label{eqdef}
     Let $\m$ be a locally finite Borel measure on $X$. Then the following holds.
     \item{$(i)$}
     $$
     \underline{\d}_{\m}(x) \leq \underline{d}_{\m}(x) \leq
     \overline{d}_{\m}(x) \leq \overline{\d}_{\m}(x).
     $$
     \item{$(ii)$} There exist the limits for $h\to\infty$ in equations 
     (\ref{eq3}), (\ref{eq4}). Moreover,
     \begin{align*}
	 \subc_{\m}(x) & = \sup_{h>0} \liminf_{t \to \infty}
	 \frac{g(t,h)} {h},  \\
	 \supc_{\m}(x) & = \inf_{h>0} \limsup_{t \to \infty}
	 \frac{g(t,h)} {h}.
     \end{align*}
     \item{$(iii)$}
     \begin{align*}
	 \subc_{\m}(x) & = \liminf_{(h,t)\to (\infty, \infty)}
	 \frac{g(t,h)} {h},  \\
	 \supc_{\m}(x) & = \limsup_{(h,t)\to (\infty, \infty)}
	 \frac{g(t,h)} {h}.
     \end{align*}
 \end{thm}

 \begin{proof}
     Properties $(i)$ and $(ii)$ follow from Proposition 1.1 in 
     \cite{GuIs9}, now we prove $(iii)$.
     Setting $ p(t,h) := \frac{f(t+h)-f(t)}{h}$, we have to 
     show that 
     $$
     \limsup_{t,h\to\infty}  p(t,h) = \lim_{h\to\infty} 
     \limsup_{t\to\infty}  p(t,h).
     $$     
     Assume $\lim_{h\to\infty}\limsup_{t\to\infty} p(t,h)=L\in\br$. 
     Let $\eps>0$, then there is $h_{\eps}>0$ such that, for any
     $h>h_{\eps}$, $\limsup_{t\to\infty} p(t,h) > L-\eps/2$, hence,
     for any $t_{0}>0$ there is $t=t(h,t_{0})>t_{0}$, such that
     $ p(t,h)>L-\eps$.  Hence, for any $h_{0}>0$, $t_{0}>0$ there
     exist $h>h_{0}$, $t>t_{0}$ such that $ p(t,h)>L-\eps$, namely
     $\limsup_{t,h\to\infty} p(t,h)\geq L$.  Conversely, assume
     $\limsup_{t,h\to\infty} p(t,h)=L'\in\br$, and choose $t_{n}$, $h_{n}$ such
     that $\lim_{n\to\infty} p(t_{n},h_{n})=L'$. 
     Then, for any 
     $h>0$, with $k_{n}$ denoting $\ceiling{\frac{h_{n}}{h}}$, we have
     \begin{align*}
	 p(t_{n},h_{n}) & \leq \frac{k_{n}h}{h_{n}}
	 p(t_{n},k_{n}h) = \frac{k_{n}h}{h_{n}}\frac{1}{k_{n}}
	\sum_{j=0}^{k_{n}-1} p(t_{n}+jh,h)\\
	&\leq \frac{k_{n}h}{h_{n}}\max_{j=0\dots k_{n}-1}
	 p(t_{n}+jh,h) \leq \frac{k_{n}h}{h_{n}}\sup_{t\geq
	t_{n}} p(t,h).
    \end{align*}
    Hence, for $n\to\infty$, we get
    $L'\leq\limsup_{t\to\infty} p(t,h)$, which implies the equality. 
    The other cases are treated analogously.
 \end{proof}

 \begin{thm}\label{subsequences}
     \begin{align*}
	 \underline{\d}_{\m}(x)
	 &=\inf_{\{t_{n}\}\to\infty}\limsup_{h\to\infty}\limsup_{n}
	 \frac{g(t_{n},h)}{h}\\
	 &=\inf_{\{t_{n}\}\to\infty}\liminf_{h\to\infty}\liminf_{n}
	 \frac{g(t_{n},h)}{h}.
     \end{align*}
     The two infima are indeed minima, indeed there exists a sequence
     $\{\underline{t}_{n}\}\to\infty$ for which
     \begin{equation*}
	 \underline{\d}_{\m}(x)
	 =\lim_{h\to\infty}\limsup_{n}
	 \frac{g(\underline{t}_{n},h)}{h}
	 =\lim_{h\to\infty}\liminf_{n}
	 \frac{g(\underline{t}_{n},h)}{h},
     \end{equation*}
     and such that any subsequence is still minimizing.
     Analogously,
     \begin{align*}
	 \overline{\d}_{\m}(x)
	 &=\sup_{\{t_{n}\}\to\infty}\limsup_{h\to\infty}\limsup_{n}
	 \frac{g(t_{n},h)}{h}\\
	 &=\sup_{\{t_{n}\}\to\infty}\liminf_{h\to\infty}\liminf_{n}
	 \frac{g(t_{n},h)}{h}.
     \end{align*}
     The two suprema are indeed maxima, indeed there exists a
     sequence $\{\overline{t}_{n}\}\to\infty$ for
     which
     \begin{equation*}
	 \overline{\d}_{\m}(x)
	 =\lim_{h\to\infty}\limsup_{n}
	 \frac{g(\overline{t}_{n},h)}{h}
	 =\lim_{h\to\infty}\liminf_{n}
	 \frac{g(\overline{t}_{n},h)}{h}.
     \end{equation*}
     and such that any subsequence is still maximizing.
 \end{thm}
 
 \begin{proof}
     We prove the second part of the Theorem, the proof of the 
     first part being analogous. Let us observe that the inequality
     $$
     \overline{\d}_{\m}(x) \geq \sup_{\{t_{n}\}\to\infty}
     \limsup_{h\to\infty} \limsup_{n} \frac{g(t_{n},h)}{h}
     $$
     obviously holds for any $t_{n}\to\infty$, therefore it is enough
     to find a sequence $\overline{t}_{n}\to\infty$ for which
     \begin{equation}\label{eqn:minor}
	 \overline{\d}_{\m}(x)\leq\liminf_{h\to\infty}\liminf_{n}
	 \frac{g(\overline{t}_{n},h)}{h}.
     \end{equation}
     In \cite{GuIs11}, Proposition 5.5, we proved that, whenever
     $dg$ is bounded by a constant $S$, where
     $$
     dg(t,h,k)=g(t, h+k) - g(t+h,k)-g(t,h),
     $$
     then, for any $\k>0$, there exists a sequence
     $\{\overline{t}_{n}\}\to\infty$ for which
     \begin{equation}\label{supinf2}
	 \limsup_{h\to+\infty} \limsup_{t\to+\infty}\frac{g(t,h)}{h}
	 \leq \liminf_{h\to+\infty}
	 \liminf_{n\in\bn}\frac{g(\overline{t}_{n},h)}{h}+\frac{2S}{\k}.
     \end{equation}
     Since in our case $dg=0$, $\{\overline{t}_{n}\}$ is the required
     sequence.  Clearly inequality (\ref{eqn:minor}) is preserved
     when passing to a subsequence.
 \end{proof}

 \subsection{Tangential dimensions  on $\br^{N}$}
 \label{subsub:TangentialFormula}
 
 In this subsection $\m$ is a Radon measure on $\br^{N}$.  Let us
 recall that the cone $T_{x}(\m)$ of measures tangent to $\m$ at
 $x\in\br^{N}$ is the set of non-zero limit points in the vague
 topology of sequences $c_{n}\m\circ D_{x}^{\l_{n}}$, where $\l_{n}$
 decreases to $0$ and $c_{n}>0$, and $D_{x}^{\l}$ is the dilation with 
 center $x\in\br^{N}$ and factor $\l>0$.  In this case vague topology is the
 weak topology determined by continuous functions with compact support.
 
 Let us consider the following properties: we say that $\m$ satisfies 
 the {\it volume doubling} condition at $x$ if
 \begin{equation}\label{voldub}
	 \limsup_{r\to0}\frac{\m(B(x,2r))}{\m(B(x,r))}<\infty,
 \end{equation}
 and the {\it weak volume doubling} condition at $x$ if there exists 
 an infinitesimal sequence $r_{n}$ such that, for any $\l>0$
 \begin{equation}\label{wvoldub}
	 \limsup_{n\to\infty}\frac{\m(B(x,\l r_{n}))}{\m(B(x,r_{n}))}<\infty.
 \end{equation}
 
 Then the following proposition holds.
 
 \begin{prop}\label{Prop:basic}
     Let $\m$ be a Radon measure on $\br^{N}$.
     \begin{itemize}
	 \item[$(i)$] $T_{x}(\mu)\neq\emptyset$ {\it iff} $\m$
	 satisfies the weak volume doubling condition at $x$.  
	 \item[$(ii)$] Assume volume doubling at $x$.  Then any
	 tangent measure of $\m$ at $x$ is of the form
	 \begin{equation} 
	     \n^{\{r_{n}\}}= \lim_{n}\frac{\m\cdot D_{x}^{r_{n}}}
	     {\m(B(x,r_{n}))},
	 \end{equation}
	 for a suitable infinitesimal sequence $\{r_{n}\}$.  Moreover, for
	 any $r_{n}\searrow0$ there exists a subsequence $r_{n_{k}}$
	 giving rise to a tangent measure $\n^{ \{ r_{n_{k}} \} }$ as above. 
	 \item[$(iii)$] Any tangent measure of $\m$ at $x$ is of the
	 form $\n^{ \{r_{n}\} }\cdot D_{x}^{\l}$ for a suitable
	 infinitesimal sequence $r_{n}$ and $\l>0$. 
	 \item[$(iv)$] Volume doubling at $x$ is equivalent to
	 $\supc_{\m}(x)<\infty$.  In particular,
	 $$
	 \limsup_{r\to0}\frac{\m(B(x,2r))}{\m(B(x,r))}=A \Rightarrow
	 \supc_{\m}(x)\leq\log_{2}A.
	 $$
     \end{itemize}
 \end{prop}
 \begin{proof}
     The first two properties are proved in \cite{Mattila}, where it is 
     also shown that any tangent measure is of the form
     $$
     \lim_{n}\frac{\m\cdot D_{x}^{r_{n}}} {\m(B(x,\l r_{n}))},
     $$
     from which $(iii)$ follows.  \\
     Let us prove $(iv)$.  By definition,
     $\limsup_{r\to0}\frac{\m(B(x,2r))}{\m(B(x,r))}=A$ can be rewritten
     as $\limsup_{t\to\infty}f(t+\log2)-f(t)=\log A$, from which
     $\supc_{\m}(x)\leq\log_{2}A$ easily follows.  Conversely, if
     $\limsup_{t\to\infty}f(t+\log2)-f(t)=\infty$, then
     $\limsup_{t\to\infty}f(t+h)-f(t)=\infty$ for any $h\geq\log2$,
     hence $\overline{\d}_{\m}(x)=\infty$.
 \end{proof}

 \begin{prop}\label{locdimtgmeas}
     Let us consider a tangent measure of $\m$ at $x$ of the form
     $\n=\n^{ \{r_{n}\} }\cdot D_{x}^{\l}$ as in Proposition
     \ref{Prop:basic}.  Then, with $f(t)=-\log\m(B(x,\e{-t}))$, and
     $t_{n}:= -\log r_{n}$, we have
     \begin{align*}
	 \underline{d}_{\n}(x)&=\liminf_{h\to\infty}\lim_{n}
	 \frac{f(t_{n}+h)-f(t_{n})}{h},\\
	 \overline{d}_{\n}(x)&=\limsup_{h\to\infty}\lim_{n}
	 \frac{f(t_{n}+h)-f(t_{n})}{h}.
     \end{align*}
 \end{prop} 
 \begin{proof}
     First let us notice that 
     $$
     \underline{d}_{\n^{ \{r_{n}\} }\cdot D_{x}^{\l}}(x)
     =\underline{d}_{\n^{ \{r_{n}\} }}(x).
     $$
     Then, let us recall that tangent measures are defined in terms of
     vague convergence, namely weak convergence on continuous
     functions with compact support.  Then, let $\f$ be a continuous
     function verifying $\chi_{B(x,1)}\leq\f\leq\chi_{B(x,2)}$. 
     Setting $h=-\log \l$, we get
     \begin{align*}
	 \n^{ \{r_{n}\} }(B(x,\l/2))
	 &=\n^{ \{r_{n}\} }\cdot  D_{x}^{\l/2}(B(x,1))\\
	 &\leq\langle\n^{ \{r_{n}\} }\cdot  D_{x}^{\l/2},\f\rangle\\
	 &=\lim_{n}\frac{\langle\m\cdot D_{x}^{\l r_{n}/2},\f\rangle} 
	 {\m(B(x,r_{n}))}\\
	 &\leq\lim_{n}\frac{\m(B(x,\l r_{n}))} 
	 {\m(B(x,r_{n}))}\\
	 &=\lim_{n}\ \exp\left[-\left(f(t_{n}+h)-f(t_{n})\right)\right]\\
	 &\leq\langle\n^{ \{r_{n}\} }\cdot  D_{x}^{\l},\f\rangle
	 \leq\n^{ \{r_{n}\} }(B(x,2\l)).
     \end{align*}
     Then, 
     \begin{equation}
	 \frac{\log\n^{ \{r_{n}\} }(B(x,\l/2))}{\log \l}
	 \leq \lim_{n}\frac{f(t_{n}+h) - f(t)}{h}
	 \leq \frac{\log\n^{ \{r_{n}\} }(B(x,2\l))}{\log \l},
     \end{equation}
     from which the thesis immediately follows.
 \end{proof}

 \begin{thm}\label{TangentialFormula}
     Let $\m$ be a Radon measure on $\br^{N}$, satisfying the volume
     doubling condition at $x$.  Then
     \begin{align*}
	 \underline{\d}_{\m}(x)&=\inf_{\n\in T_{x}(\m)}
	 \underline{d}_{\n}(x), \\
	 \overline{\d}_{\m}(x)&=\sup_{\n\in T_{x}(\m)} 
	 \overline{d}_{\n}(x).
     \end{align*}
 \end{thm}
 
 Let us remark that volume doubling implies weak volume doubling, 
 namely the set of tangent measures at $x$ is non-empty.
 
 \begin{proof}
     Let us give the proof for $\underline{\d}_{\m}(x)$, the other 
     case being proved analogously. Let $T(\m,x)$ be the set of sequences 
     $t_{n}\to\infty$ such that $r_{n}=e^{-t_{n}}$ generates a 
     tangent measure as in Lemma \ref{Prop:basic} $(ii)$. Then, from 
     Proposition \ref{locdimtgmeas}, we get
     \begin{equation}
	 \inf_{\n\in T_{x}(\m)} \underline{d}_{\n}(x)
	 =\inf_{\{t_{n}\}\in T(\m,x)}\liminf_{h\to\infty}\liminf_{n}
	 \frac{f(t_{n}+h)-f(t_{n})}{h}.
     \end{equation}
     So the equality is proved if we show that $T(\m,x)$ contains one
     of the minimizing sequences of theorem \ref{subsequences}.  This
     holds true, since any minimizing sequence of theorem
     \ref{subsequences} has a subsequence giving rise to a tangent
     measure by Proposition \ref{Prop:basic} $(ii)$, and such
     subsequence inherits the minimizing property.
 \end{proof}

 \subsection{Further properties}
 Tangential dimensions are invariant under bi-Lipschitz maps.
 
 \begin{prop}
     Let $X,\ Y$ be metric spaces, $f:X \to Y$ be a bi-Lipschitz map
     $i.e.$ there is $L>0$ such that $L^{-1} d_{X}(x,x') \leq
     d_{Y}(f(x),f(x')) \leq L d_{X}(x,x')$, for $x,x'\in X$.  Let $\m$
     be a finite Borel measure on $X$ and set $\n:= \m\circ f^{-1}$,
     which is a finite Borel measure on $Y$.  Then
     $\underline{\d}_{\m}(x) = \underline{\d}_{\n}(f(x))$ and
     $\overline{\d}_{\m}(x) = \overline{\d}_{\n}(f(x))$, for all $x\in
     X$.
 \end{prop}
 \begin{proof}
     Observe that, for any $x\in X$, $y\in Y$, $r>0$, we have 
     \begin{align*}
	 B(f(x),r/L) & \subset f(B(x,r)) \subset B(f(x),rL) \\
	 B(f^{-1}(y),r/L) & \subset f^{-1}(B(y,r)) \subset
	 B(f^{-1}(y),rL)
     \end{align*}
     which implies
     $$
     \m(B(x,R/L)) \leq \n(B(f(x),R)) = \m(f^{-1}(B(f(x),R)))
     \leq \m(B(x,RL)).
     $$
     Therefore, taking $\limsup_{R\to 0}$, then $\lim_{\l\to0}$, and
     doing some algebra, we get
     $$
     \lim_{\l\to0} \limsup_{R\to 0} \frac{\log
     \frac{\m(B(x,R))}{\m(B(x,\l R))}}{\log L^{2}/\l} \leq
     \lim_{\l\to0} \limsup_{R\to 0} \frac{\log
     \frac{\n(B(x,R))}{\n(B(x,\l R))} }{\log 1/\l} \leq \lim_{\l\to0}
     \limsup_{R\to 0} \frac{\log \frac{\m(B(x,R))}{\m(B(x,\l R))}}{\log
     1/(L^{2}\l)}
     $$
     which means $\overline{\d}_{\m}(x) = \overline{\d}_{\n}(f(x))$. 
     The other equality is proved in the same manner.
 \end{proof}
 
 The following propositions show some properties of tangential 
 dimensions, $i.e.$ their behaviour under the operations of sum or 
 tensor product of measures. 
 
 \begin{prop}
     Let $\m_{1},\ \m_{2}$ be finite Borel measures on $X$.  Then
     \begin{align*}
	\underline{\d}_{\m_{1}+\m_{2}}(x) & \geq \min \{
	\underline{\d}_{\m_{1}}(x) , \underline{\d}_{\m_{2}}(x) \} \\
	\overline{\d}_{\m_{1}+\m_{2}}(x) & \leq \max \{
	\overline{\d}_{\m_{1}}(x) , \overline{\d}_{\m_{2}}(x) \}.
    \end{align*}
 \end{prop}
 \begin{proof}
     As
     $$
     \frac{1}{2}\ \min \left\{ \frac{a}{c}, \frac{b}{d} \right\} \leq
     \frac{a+b}{c+d} \leq 2\ \max \left\{ \frac{a}{c}, \frac{b}{d}
     \right\}
     $$
     we get
     \begin{align*}
	-\frac{\log 2}{\log 1/\l} + \ & \min \left\{
	\frac{\log\frac{\m_{1}(B(x,r))}{\m_{1}(B(x,\l r))}}{\log 1/\l} ,
	\frac{\log \frac{\m_{2}(B(x,r))}{\m_{2}(B(x,\l r))}}{\log 1/\l}
	\right\} \\
	& \leq \frac{\log
	\frac{\m_{1}(B(x,r))+\m_{2}(B(x,r))}{\m_{1}(B(x,\l
	r))+\m_{2}(B(x,\l r))}}{\log 1/\l} \\
	\leq \frac{\log 2}{\log 1/\l} + \ & \max \left\{
	\frac{\log\frac{\m_{1}(B(x,r))}{\m_{1}(B(x,\l r))}}{\log 1/\l} ,
	\frac{\log \frac{\m_{2}(B(x,r))}{\m_{2}(B(x,\l r))}}{\log 1/\l}
	\right\}.
    \end{align*}
    Therefore, taking $\limsup_{r\to0}$, then using the equality
    $\limsup_{r\to0} \max \{ f(r), g(r) \} = \max \{ \limsup_{r\to0}
    f(r) , \limsup_{r\to0} g(r) \} $, and finally taking
    $\lim_{\l\to0}$, we obtain
    $$
    \overline{\d}_{\m_{1}+\m_{2}}(x) \leq \max \{
    \overline{\d}_{\m_{1}}(x) , \overline{\d}_{\m_{2}}(x) \}.
    $$
    Besides, taking $\liminf_{r\to0}$, then using the equality $
    \liminf_{r\to0} \min \{ f(r), g(r) \} = \min \{ \liminf_{r\to0}
    f(r) , \liminf_{r\to0} g(r) \} $, and finally taking
    $\lim_{\l\to0}$, we obtain
    $$
    \underline{\d}_{\m_{1}+\m_{2}}(x) \geq \min \{
    \underline{\d}_{\m_{1}}(x) , \underline{\d}_{\m_{2}}(x) \}.
    $$
 \end{proof}

 \begin{prop}
     Let $X,\ Y$ be metric spaces, $\m,\ \n$ finite Borel measures on 
     $X$ and $Y$ respectively. Then 
     \begin{align*}
         \underline{\d}_{\m\otimes\n}((x,y)) & \geq
         \underline{\d}_{\m}(x) + \underline{\d}_{\n}(y)  \\
         \overline{\d}_{\m\otimes\n}((x,y)) & \leq \overline{\d}_{\m}(x) + 
         \overline{\d}_{\n}(y).
     \end{align*}
 \end{prop}
 \begin{proof}
     Endow $X\times Y$ with the metric
     \begin{equation}\label{eqn:prod.metric}
	 d((x_{1},y_{1}),(x_{2},y_{2})) := \max \{ d_{X}(x_{1},x_{2}),
	 d_{Y}(y_{1},y_{2}) \}
     \end{equation}
     which is by-Lipschitz equivalent to the product metric.  Then
     \begin{equation}\label{eqn:prod.ball}
	 B_{X\times Y}((x,y),R) = B_{X}(x,R) \times B_{Y}(y,R),
     \end{equation}
     which implies
     $$
     \m\otimes\n(B_{X\times Y}((x,y),R)) = \m(B_{X}(x,R)) \n(B_{Y}(y,R)),
     $$
     and
     $$
     \frac{\log\frac{\m\otimes\n(B_{X\times
     Y}((x,y),r))}{\m\otimes\n(B_{X\times Y}((x,y),\l r))}}{\log 1/\l} =
     \frac{\log\frac{\m(B_{X}(x,r))}{\m(B_{X}(x,\l r))}}{\log 1/\l} +
     \frac{\log\frac{\n(B_{ Y}(y,r))}{\n(B_{Y}(y,\l r))}}{\log 1/\l},
     $$
     from which the thesis follows.
 \end{proof}

 The following theorem examines the dependence of tangential dimensions 
 on the point $x\in X$.
 
 \begin{thm}\label{thm:misurabile}
     The function 
     $$
     \underline{\d}_{\m}: x\in X \to \lim_{\l\to0}
     \liminf_{r\to0}\frac{\log\left( \frac{\m(B(x,r))}{\m(B(x,\l r))}
     \right)}{\log 1/\l} \in [0,\infty)
     $$ 
     is Borel-measurable.  The same is true of $\overline{\d}_{\m}$
     with $\liminf$ replaced by $\limsup$.
 \end{thm}
 \begin{proof}
     Set, for $r>0$, $\l\in(0,1)$, $f_{r,\l}(x) :=
     \frac{\m(B(x,r))}{\m(B(x,\l r))}$, which is Borel-measurable by
     \cite{Edgar}, proof of 1.5.9.  Then we must prove that
     $$
     f(x) := \lim_{\l\to0} \frac{1}{\log 1/\l}
     \log\left(\liminf_{r\to0} f_{r,\l}(x) \right),
     $$
     is Borel-measurable.  First 
     $$
     f_{\l}(x) := \liminf_{r\to0} f_{r,\l}(x) =  \lim_{r\to 0} \inf_{0<r'<r}
     f_{r',\l}(x)
     $$ 
     is Borel-measurable, because, from $\{r_{n}\}\subset\bq$,
     $r_{n}\nearrow r$, it follows $f_{r_{n},\l}(x) \to f_{r,\l}(x)$,
     and
     $$
     \lim_{r\to 0} \inf_{0<r'<r} f_{r',\l}(x) = \lim_{r\to 0}
     \inf_{0<r'<r, r'\in\bq} f_{r',\l}(x) = \lim_{n \to \infty}
     \inf_{0<r'<\frac{1}{n}} f_{r',\l}(x).
     $$
     Then 
     $$
     f(x) = \lim_{\l\to0} \frac{\log f_{\l}(x)}{\log 1/\l} =
     \lim_{n\to \infty} \frac{\log f_{\frac{1}{n}}(x)}{\log n}
     $$ 
     is Borel-measurable.
 \end{proof}

\subsection{Relations between tangential dimensions of metric spaces and measures}

 \begin{dfn}
     Let $(X,d)$ be a metric space, $E\subset X$, $x\in E$.  Let us
     denote by $n(r,E) \equiv n_{r}(E)$ the minimum number of open
     balls of radius $r$ necessary to cover $E$, and by $\n(r,E)
     \equiv \n_{r}(E)$ the maximum number of disjoint open balls of
     $E$ of radius $r$ contained in $E$. We call {\it upper, resp.  lower
     tangential dimension} of $E$ at $x$ the (possibly infinite)
     numbers
     \begin{align*}
	\underline{\d}_{E}(x) & := \liminf_{\l \to 0} \liminf_{r \to
	0} \frac{\log n(\l r, E\cap\B(x,r))}{\log 1/\l}, \\
	\overline{\d}_{E}(x) & := \limsup_{\l \to 0} \limsup_{r \to
	0} \frac{\log n(\l r, E\cap\B(x,r))}{\log 1/\l}.
    \end{align*}
 \end{dfn}
 Observe that we obtain the same numbers if we use $\n$ in place of 
 $n$.

\begin{thm}[\cite{GuIs11}]\label{newformula}
	Let us assume the following conditions 
	\itm{i}
	$\displaystyle
	\limsup_{r\to0}n_{\l
	r}(\ov{B}_{X}(x,r))<\infty\quad \forall \l>0,$
    \itm{ii}
	there exist constants $c\geq1$, $a\in(0,1]$ such that, for any
	$r\leq a$, $\l,\m\leq1$, $y,z\in B_{X}(x,r)$,
	$\displaystyle
	    n(\l\m r, B_{X}(y,\l r))\leq c n(\l\m r, B_{X}(z,\l r)).$
	    \\
    \noindent Then
    \begin{align}
	\subc_{X}(x) &= \inf_{T\in\ct_{x} X}\subd(T)=
	\inf_{T\in\ct_{x} X}\supd(T),\label{eq.5}\\
	\supc_{X}(x) &= \sup_{T\in\ct_{x} X}\subd(T)=
	\sup_{T\in\ct_{x} X}\supd(T).\label{eq.6}
    \end{align}
\end{thm}

\begin{lem}
    Let $\m$ be a finite Borel measure on the metric space $X$, $x\in
    X$.  The following inequalities hold:
    \begin{align}
	\m(B(x,r))
	&\leq n(\l r,B(x,r)) 
	\sup_{y\in B(x,r)}\m(B(y,\l r)).\label{n-mu-ineq}\\
	\m(B(x,r))
	&\geq \n(\l r,B(x,r)) 
	\inf_{y\in B(x,r)}\m(B(y,\l r)).\label{ni-mu-ineq}
    \end{align}
\end{lem}

\begin{proof}
    The first inequality follows from 
    $$
    B(x,r) \subset \bigcup_{i=1}^{n(\l r,B(x,r))} B(y_{i},\l r),
    $$
    the second follows from
    $$
    B(x,r) \supset \bigcup_{i=1}^{\n(\l r,B(x,r))} B(y_{i},\l r).
    $$
\end{proof}

 \begin{prop}\label{prop:ptw.box=ptw.mis}
     Let $\m$ be a finite Borel measure on $X$ and define
     \begin{align*}
	 m_{x}(r,R) & := \inf \{ \m(B(y,r) : B(y,r) \subset B(x,R) \}
	 \\
	 M_{x}(r,R) & := \sup \{ \m(B(y,r) : B(y,r) \subset B(x,R) \}.
     \end{align*}
     If 
     $$
     \limsup_{\l \to 0} \limsup_{r \to 0} \frac{ \log \frac{ M_{x}(\l
     r,r) }{ m_{x}(\l r,r) }}{ \log 1/\l} = 0
     $$
     then $\underline{\d}_{X}(x) = \underline{\d}_{\m}(x)$ and
     $\overline{\d}_{X}(x) = \overline{\d}_{\m}(x)$.
 \end{prop}
 \begin{proof}
     From the definition of $M_{x}$ we get
     $$
     \m(B(x,R)) \leq \sum_{i=1}^{n_{r}(B(x,R))} \m(B(y_{i},r)) \leq
     n_{r}(B(x,R)) M_{x}(r,R+r).
     $$
     From the definition of $m_{x}$ we get
     $$
     \m(B(x,R)) \geq \sum_{i=1}^{\n_{r}(B(x,R))} \m(B(y_{i},r)) \geq
     \n_{r}(B(x,R)) m_{x}(r,R).
     $$
     Therefore
     $$
     \n_{\l R}(B(x,R))\ \frac{ m_{x}(\l R,R) }{ M_{x}(\l R,R) } \leq
     \frac{\m(B(x,R))}{\m(B(x,\l R))} \leq n_{\l R}(B(x,R))\ \frac{
     M_{x}(\l R,(\l + 1)R) }{ m_{x}(\l R,(\l + 1)R) }
     $$
     and
     \begin{align*}
	 \frac{\log \n_{\l R}(B(x,R))}{\log 1/\l} - \frac{\log \frac{
	 M_{x}(\l R,R) }{ m_{x}(\l R,R) }}{\log 1/\l} & \leq
	 \frac{\log \frac{\m(B(x,R))}{\m(B(x,\l R))}}{ \log 1/\l} \\
	 &\leq \frac{\log n_{\l R}(B(x,R))}{\log 1/\l} + \frac{\log
	 \frac{ M_{x}(\l R,(\l + 1)R) }{ m_{x}(\l R,(\l + 1)R) }}{\log
	 1/\l}.
     \end{align*}
     Taking $\limsup_{R \to 0}$ and then $\limsup_{\l \to 0}$ and
     doing some algebra we get
     \begin{align*}
	& \limsup_{\l \to 0} \limsup_{R \to 0} \frac{\log\n_{\l
	R}(B(x,R))}{\log 1/\l} - \liminf_{\l \to 0} \liminf_{R \to 0}
	\frac{\log \frac{ M_{x}(\l R,R) }{ m_{x}(\l R,R) }}{\log 1/\l}
	\\
	& \leq \lim_{\l \to 0} \limsup_{R \to 0}  \frac{\log
	\frac{\m(B(x,R))}{\m(B(x,\l R))}}{ \log 1/\l} \\
	& \leq \limsup_{\l \to 0} \limsup_{R \to 0} \frac{\log n_{\l
	R}(B(x,R))}{\log 1/\l} + \limsup_{\l \to 0} \limsup_{R \to 0}
	\frac{\log \frac{ M_{x}(\l R,R) }{ m_{x}(\l R,R) }}{\log 1/\l}
    \end{align*}
    and the thesis $\overline{\d}_{X}(x) = \overline{\d}_{\m}(x)$
    follows.  The proof of the other equality is analogous.
 \end{proof}
 
 \begin{cond}\label{NewAssum}
	 Let $\m$ be a finite Borel measure on $X$. For any $x\in X$ there
	 are constants $R,\ C>0$, depending on $x$, such that, for any
	 $y\in B(x,R)$, $r\in (0,R)$, it holds
	 $$
	 C^{-1} \m(B(x,r)) \leq \m(B(y,r)) \leq C \m(B(x,r)).
	 $$
 \end{cond}
  
 \begin{cor}\label{cor:ptw.box=ptw.mis}
	 Let $\m$ be a finite Borel measure on $X$ satisfying
	 Condition \ref{NewAssum}.  Then $\underline{\d}_{X}(x) =
	 \underline{\d}_{\m}(x)$ and $\overline{\d}_{X}(x) =
	 \overline{\d}_{\m}(x)$, and these functions are locally
	 constant.
 \end{cor}
 \begin{proof}
	 As $m_{x}(\l r,r) = \inf \{ \m(B(y,\l r) : B(y,\l r)
	 \subset B(x,r) \} \geq C^{-1} \m(B(x,\l r))$, and
	 $M_{x}(\l r,r) \leq C \m(B(x,\l r))$, for any 
	 $\l\in(0,1)$, $r\in(0,R)$, we get
	 $$
	 \limsup_{\l \to 0} \limsup_{r \to 0} \frac{ \log \frac{ 
	 M_{x}(\l r,r) }{ m_{x}(\l r,r) }}{ \log 
	 1/\l} \leq \limsup_{\l\to0} \frac{C^{2}}{\log 1/\l} = 0,
	 $$
	 and the thesis follows from the previous Proposition. Moreover, 
	 $\forall y\in B(x,R)$, 
	 $\underline{\d}_{\m}(y) = \underline{\d}_{\m}(x)$, and 
	 $\overline{\d}_{\m}(y) = \overline{\d}_{\m}(x)$.
 \end{proof}

 We now show that Condition \ref{NewAssum} also implies property
 $(ii)$ of Theorem \ref{newformula} and volume doubling, hence
 tangential dimensions are indeed suprema, resp.  infima, of
 dimensions of tangent objects.  We first need some Lemmas.
 
 \begin{lem}\label{number-ineq}
     The following inequality holds, for $0\leq \l,\m \leq1$:
     \begin{equation*}
	 n(\l\m r,B_{X}(x,r))
	 \leq n(\l r,B_{X}(x,r)) 
	 \sup_{y\in B_{X}(x,r)}n(\l\m r,B_{X}(y,\l r)).
     \end{equation*}
 \end{lem}

 \begin{proof}
     Let us note that we may realize a covering of $B_{X}(x,r)$
     with balls of radius $\l\m r$ as follows: first choose an
     optimal covering of $B_{X}(x,r)$ with balls of radius
     $\l r$, and then cover any covering ball optimally with balls
     of radius $\l\m r$.  The thesis follows.    
 \end{proof}
 
 \begin{lem}\label{fin-dim}
     Let $X$ be a subset of $\br^{N}$. Then, for any $\l\leq1$, there 
     exists a constant $K_{\l}$ such that
     $$
     n(\l r,B_{X}(x,r))\leq K_{\l},\qquad \forall r>0, x\in X.
     $$
 \end{lem}    
 \begin{proof}
     Since (cf.  e.g. \cite{GuIs6}) the inequality
     \begin{equation}\label{n-ni}
	 n(r,B_{X}(x,R)) \geq
	 \n(r,B_{X}(x,R)) \geq n(2r,B_{X}(x,R))
     \end{equation}
     holds, we get 
     \begin{equation*}
	 n(\l r,B_{X}(x,r))\leq \n\left(\frac{\l}{2}
	 r,B_{X}(x,r)\right) \leq \n\left(\frac{\l}{2}
	 r,B_{\br^{N}}(x,r)\right)=\n\left(\frac{\l}{2},B_{\br^{N}}(1)\right),
     \end{equation*}
     where we used the dilation invariance of $\br^{N}$ in the last
     equation, and omitted the irrelevant reference to the point $x$ in
     the last term.
 \end{proof}
 
 \begin{lem}                                                      
     Let $X$ be a closed subset of $\br^{N}$, $x\in X$.  Then                
     property $(ii)$ of Theorem \ref{newformula} is equivalent to the following:                
     \\                                                                      
     For any $\th>0$, or equivalently for some $\th>0$, there exist
     constants $c_{\th}\geq1$, $a_{\th}\in(0,1]$ such that, for any
     $r\leq a_{\th}$, $\l,\m\leq1$, $y,z\in B_{X}(x,r)$,
     \begin{equation}\label{eq:ass-theta}                                    
	 n(\l\m r, B_{X}(y,\l r))\leq c_{\th} n(\l\m r\th, B_{X}(z,\l r)).   
     \end{equation}                                                          
 \end{lem}                                                                 
 
 \begin{proof}                                                               
     First we show that property $(ii)$ of Theorem \ref{newformula}
     implies (\ref{eq:ass-theta}) for $\th>1$, hence for all $\th>0$. 
     Indeed, for $a_{\th}=\frac{a}{\th}$, we get
     \begin{align*}                                                          
	 \frac{n(\l\m r, B_{X}(y,\l r))}{n(\l\m r\th, B_{X}(z,\l r))}        
	 &\leq c^{2}                                                         
	 \frac{n(\l\m r, B_{X}(x,\l r))}{n(\l\m r\th, B_{X}(x,\l r))}\\      
	 &\leq c^{3} n(\l\m r, B_{X}(x,\l \m r\th))                          
	 \leq c^{3}K_{1/\th},                                                
     \end{align*}                                                            
     where we used Lemmas \ref{number-ineq}, \ref{fin-dim}.                  
     We get (\ref{eq:ass-theta}) with $c_{\th}=c^{3} K_{1/\th}$.             
     \\                                                                      
     Now we prove (\ref{eq:ass-theta}), for some $\th<1$, implies
     property $(ii)$ of Theorem \ref{newformula}.  We set $a=\th
     a_{\th}$.  Then, reasoning as in the previous case, we get
     \begin{align*}                                                          
	 \frac{n(\l\m r, B_{X}(y,\l r))}{n(\l\m r, B_{X}(z,\l r))}           
	 &\leq c_{\th}^{3} n(\l\m r\th^{2}, B_{X}(x,\l \m r/\th))\\              
	 &\leq c_{\th}^{3} K_{\th^{3}}.                                      
     \end{align*}
     Finally we observe that (\ref{eq:ass-theta}), for some $\th>0$,
     implies (\ref{eq:ass-theta}), for some $\th<1$, hence, because of
     what has already been proved, it implies (\ref{eq:ass-theta}),
     for all $\th>0$.  The thesis follows.
 \end{proof}                                                                 
                                                                                
 \begin{prop}\label{MuImplyN}
     Let $\m$ be a finite Borel measure on $X$ satisfying Condition
     \ref{NewAssum}.  Then property $(ii)$ of Theorem \ref{newformula}
     and Volume Doubling hold.
 \end{prop}
 
 \begin{proof}
     Let us show volume doubling. Indeed, if $r\leq R/2$,
     $$
     \m(B(x,2r))\leq n(r,B(x,2r))\sup_{y\in B(x,2r)}\m(B(y,r))
     \leq C n(r,B(x,2r)) \m(B(x,r)),
     $$
     hence, by Lemma \ref{fin-dim},
     $$
     \frac{\m(B(x,2r))}{\m(B(x,r))}\leq C  n(r,B(x,2r))
     \leq C K_{1/2}.
     $$
     Now we prove property $(ii)$ of Theorem \ref{newformula} in the
     equivalent form (\ref{eq:ass-theta}), for $\th=1/2$, with
     $c_{1/2}=C^{4}$ and $a_{1/2}=R/2$.  Indeed, let $\l,\m\leq1$,
     $r\leq R/2$, $y,z\in B(x,r)$.  By (\ref{n-mu-ineq}), we get
     \begin{align*}
	 n(\l\m r/2,B(z,\l r))
	 &\geq \frac{\m(B(z,\l r))}
	 {\displaystyle{\sup_{w\in B(z,\l r)}}\m(B(w,\l\m r/2))}\\
	 &\geq \frac{\m(B(z,\l r))}
	 {\displaystyle{\sup_{w\in B(z, 2 r)}}\m(B(w,\l\m r/2))}
	 \geq\frac{1}{C^{2}}\frac{\m(B(x,\l r))}{\m(B(x,\l\m r/2))}.
     \end{align*}
     Analogously, by (\ref{ni-mu-ineq}), we get
     \begin{align*}
	 \n(\l\m r/2,B(y,\l r))
	 &\leq \frac{\m(B(y,\l r))}
	 {\displaystyle{\inf_{w\in B(y,\l r)}}\m(B(w,\l\m r/2))}\\
	 &\leq \frac{\m(B(y,\l r))}
	 {\displaystyle{\inf_{w\in B(y, 2 r)}}\m(B(w,\l\m r/2))}
	 \leq C^{2}\frac{\m(B(x,\l r))}{\m(B(x,\l\m r/2))}.
     \end{align*}
     Finally, making use of (\ref{n-ni}), we get
     \begin{equation*}
	 \frac{n(\l\m r,B(y,\l r))}{n(\l\m r/2,B(z,\l r))}\leq C^{4}.
     \end{equation*}     
 \end{proof}

 \begin{cor}\label{cor:new}
     If Condition \ref{NewAssum} holds for a measure $\m$ on $F$, then
     tangential dimensions for $\m$ are extrema of local dimensions of
     tangent measures, and tangential dimensions for $F$ are extrema
     of local dimensions of tangent sets.
 \end{cor}
  
 \section{Computation of tangential dimensions}
 \subsection{Self-similar fractals}
 We compute here the tangential dimensions for self-similar fractals
 with open set condition, showing that Condition \ref{NewAssum} is
 satisfied, and that upper and lower tangential dimensions for the
 Hausdorff measure are equal, hence coincide with the Hausdorff
 dimension. This means that self-similar fractals are too regular to 
 give rise to a dimension interval, and a different class has to be 
 considered, see the next subsection.
 
 Let us recall that a self-similar fractal $F$ is the fixed point of a 
 map
 \begin{align*}
     W : K\in\cpt{\br^{N}} \to \cup_{j=1}^{p} w_{j}(K) \in    
     \cpt{\br^{N}},
 \end{align*}
 where $w_{j}$ are similarities with similarity parameter $\l_{j}$,
 and that it satisfies the {\it open set condition} if there exists an
 open set $V$ s.t. $w_{j}V\subset V$.  It is well known that the
 Hausdorff dimension $d$ of $F$ satisfies 
 $\sum_{j=1}^{p}\l_{j}^{d}=1$, that the corresponding Hausdorff 
 measure $\ch_{d}$ is nontrivial on $F$ and the normalized restriction 
 of $\ch_{d}$ to $F$ is the unique self-similar probability measure on 
 $F$. In particular $\ch_{d}(w_{j}(F))=\l_{j}^{d}$.
 
 If $\S_{n}$ is the set of multi-indices $\s=(\s_{1}, \dots, \s_{n})$
 of length $n$, we denote by $w_{\s}$ the product
 $w_{\s_{1}}\cdot\dots\cdot w_{\s_{n}}$, and by $\l_{\s}$ the product
 $\l_{\s_{1}}\cdot\dots\cdot\l_{\s_{n}}$.  Also we use the notations
 $F_{\s}=w_{\s}F$, $V_{\s}=w_{\s}V$.  We note that
 $\ch_{d}(F_{\s})=\ch_{d}(V_{\s})=\l_{\s}^{d}$.
 
 \begin{lem}\label{sscond}
     Let $F$ be a self-similar fractal with open set condition. There 
     exists a constant $C>0$ s.t., for any $x\in F$, $r>0$,
     $$
     C^{-1} r^{d}\leq \ch_{d}(B_{F}(x,r))\leq C r^{d}.
     $$
 \end{lem}
 
 \begin{proof}
     It is not restrictive to assume that the diameter of $V$ is $1$,
     hence $V_{\s}$ has diameter $\l_{\s}$.  If $\s$ is a multi-index
     of length $n$, in the following we shall denote by $\ov{\s}$ the
     multi-index $(\s_{1},\dots,\s_{n-1})$.  Let us consider the set
     $\S(r)$ of multi-indices $\s$ s.t. $\l_{\s}<r\leq \l_{\ov{\s}}$.
     Clearly the $F_{\s}$'s, $\s\in\S(r)$ give a covering of $F$, and
     the $V_{\s}$, $\s\in\S(r)$, are pairwise disjoint.  Then, if
     $x\in F_{\s}$, $\s\in\S(r)$, $B_{F}(x,r)\supset F_{\s}$,
     whence
     $$
     \ch_{d}(B_{F}(x,r))\geq \l_{\s}^{d}\geq r^{d} \underline{\l}^{d},
     $$
     where $\underline{\l}=\min(\l_{1},\dots\l_{p})$.
     Conversely, set $\S(r,x)=\{\s\in\S(r):F_{\s}\cap 
     B_{F}(x,r)\ne\emptyset\}$. Then
     $$
     \ch_{d}(B_{F}(x,r))\leq \sum_{\s\in\S(r,x)}\l_{\s}^{d}
     \leq \#\S(r,x)\,r^{d}.
     $$
     Observe now that $\cup_{\s\in\S(r,x)}V_{\s}\subset B(x,2r)$, 
     therefore
     $$
     \omega_{N}(2r)^{N}:=\vol(B(x,2r))\geq \sum_{\s\in\S(r,x)} 
     \vol(V)\l_{\s}^{N}\geq \vol(V)\,\#\S(r,x)\,r^{N} \underline{\l}^{N},
     $$
     therefore
     $$
     \#\S(r,x) \leq \frac{\omega_{N}2^{N}}{\vol(V)\,
     \underline{\l}^{N}}.
     $$
     The thesis follows.
 \end{proof}
 
 \begin{cor}
     Let $F$ be a self-similar fractal satisfying open set condition,
     $d$ its Hausdorff dimension.  Then, for any $x\in F$,
     $$
     d=\subc_{F}(x)=\supc_{F}(x)=\subc_{\ch_{d}}(x)=\supc_{\ch_{d}}(x).
     $$     
 \end{cor}
 
 \begin{proof}
     Since Condition \ref{NewAssum} is satisfied for $\ch_{d}$, it is
     enough to compute the tangential dimensions relative to the
     Hausdorff measure.  Indeed Lemma \ref{sscond} may be rephrased as
     $d\,t-\log C\leq f(t)\leq d\,t+\log C$, where $f(t) = -\log 
     \ch_{d}(B(x,\e{-t}))$, from which it follows
     $$
     \lim_{h\to\infty}\lim_{t\to\infty}\frac{f(t+h)-f(t)}{h}=d.
     $$
 \end{proof}

 \subsection{Translation fractals}
 Now we compute the tangential dimensions for translation fractals
 defined in \cite{GuIs8} (cf.  also \cite{GuIs9,GuIs10}).
 
 Let $\{w_{nj} \}$, $n\in\bn$, $j=1,\ldots,p_{n}$, be contracting
 similarities of $\br^{N}$, with contraction parameter
 $\l_{n}\in(0,1)$ only depending on $n$, and assume they verify the
 {\it regular open set condition}, namely there exists a nonempty
 bounded open set $V$ in $\br^N$ for which $w_{nj}(V) \subset V$, the 
 Lebesgue measure of $V$ is equal to the Lebesgue measure of its 
 closure $C$ and $V$ is equal to the interior of $C$. 
 Setting $W_{n} : K\in\cpt{\br^{N}} \to \cup_{j=1}^{p_{n}} w_{nj}(K)
 \in \cpt{\br^{N}}$, we get a sequence of compact sets $\{ W_{1}\circ
 W_{2}\circ\cdots\circ W_{n}(\ov{V}) \}$ contained in $\ov{V}$, we
 call the Hausdorff limit $F$ a {\it translation fractal}.  Since the
 sequence is indeed decreasing, $F$ can be equivalently defined as the
 intersection.  To avoid triviality we assume $p_{j}\geq2$, which
 implies $2\l_{j}^{N}\leq1$, i.e. $\l_{j}\leq 2^{-1/N}$.
 
 As an example of fractals in our class consider the following 
 construction. It is obtained by applying a sequence of either a 
 Carpet step or a Vicsek step.

 The Carpet step ($q=1$) is obtained by dividing the sides of a square in $3$
 equal parts, so as to obtain $9$ equal squares, and then the central square
 is removed.

 The Vicsek step ($q=2$) is obtained by dividing the sides
 of a square in $3$ equal parts, so as to obtain $9$
 equal squares, and then $4$ squares are removed, so that to remain
 with a chessboard.

 In particular, we may set $q_{j}=1$, if $(k-1)(2k-1)<j\leq (2k-1)k$
 and $q_{j}=2$, if $ k(2k-1)<j\leq k(2k+1)$, $k=1,2,\dots$, getting a
 translation fractal with dimensions given by (cf.  Theorem
 \ref{TgDimForLimFrac} below)
 $$
 \subc=\frac{\log5}{\log3}<\subd=\supd=\frac{\log40}{\log9}
 <\supc=\frac{\log8}{\log3}
 $$

\begin{figure}[ht]
     \centering
     \subfigure{
     \psfig{file=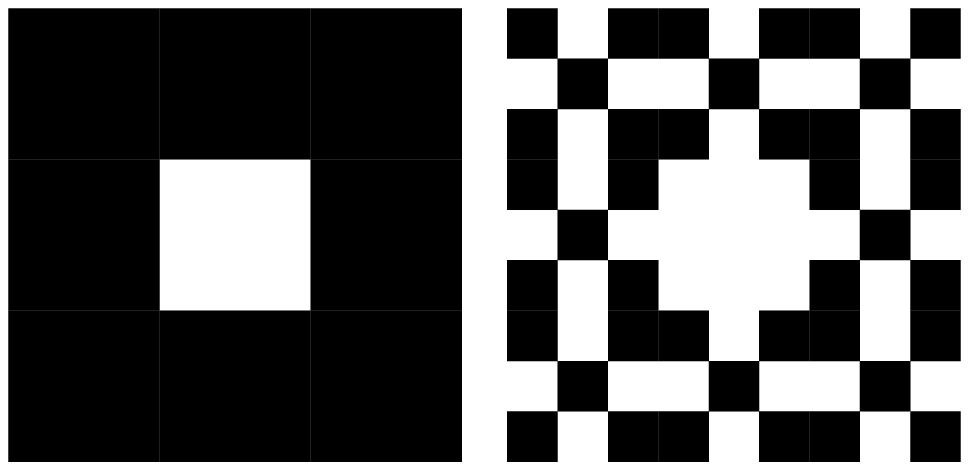,height=1.5in}}
     \hspace{0.3 in}
     \subfigure{
     \psfig{file=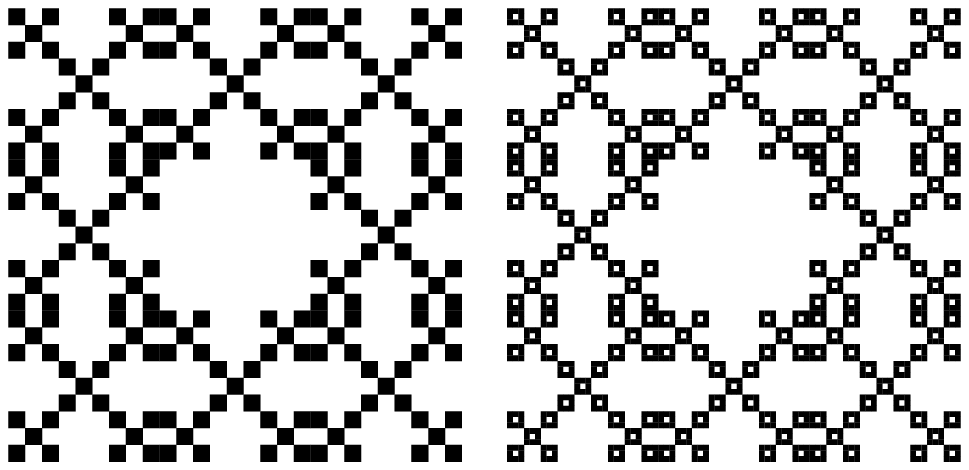,height=1.5in}}
     \caption{Carpet-Vicsek}
     \label{fig:Carpet-Vicsek}
 \end{figure}

 The first four steps ($q=1,2,2,1$) of the procedure above are shown
 in Figure \ref{fig:Carpet-Vicsek}.
 
 More examples are contained in \cite{GuIs11}.

 We set $\La_{n}=\prod_{i=1}^{n}\l_{i}$, $P_{n}=\prod_{i=1}^{n}p_{i}$,
 $\Sigma_{n} := \{\s:\{1,\ldots,n\}\to \bn : \s(k) \in
 \{1,\ldots,p_{k}\}, k=1,\ldots,n \}$, $\Sigma := \cup_{n\in\bn}
 \Sigma_{n}$, and write $w_{\s} := w_{1\s(1)}\circ w_{2\s(2)} \circ
 \cdots \circ w_{n\s(n)}$, for any $\s\in\Sigma_{n}$, and
 $V_{\s}:=w_{\s}V$, $C_{\s}:=w_{\s}C$.
 
 On $F$ there is a canonical limit measure $\m$, that can be defined
 as the weak$^{*}$-limit of the sequence
 \begin{equation}
     \m_{n}(A) =\sum_{|\s|=n}P_{n}^{-1}\m_{0}(w_{\s}^{-1}(A)),
 \end{equation}
 the limit being independent of the starting probability Borel measure
 $\m_{0}$.  In \cite{GuIs10}, Theorem 1.7, we proved that, when $V$ is
 regular, $\m$ can be characterized via the following property: for any
 subset $\cai$ of $\S_{n}$
 \begin{equation}\label{ineqCsa}
     \m(V_{\cai}) \leq \frac{\#\cai}{P_{n}} \leq
     \m(C_{\cai}),
 \end{equation}
 where we set $C_{\cai}=\cup_{\s\in\cai}C_{\s}$, $V_{\cai}$ equal 
 to the interior of $C_{\cai}$ relative to $C$.
 
 \begin{thm}\label{ExtrTgDim}
     Let $F$ be a translation fractal with regular open set condition, with
     the notation above, and assume $p:=\sup_{n}p_{n}<\infty$.  Then
     Condition \ref{NewAssum} holds for the limit measure $\m$,
     therefore tangential dimensions for $F$ and for $\m$ coincide, 
     and they are extremal  dimensions for the corresponding tangent objects.
 \end{thm} 

 \begin{proof}
     Denote by $B_{F}(x,r) := F\cap B(x,r)$, $x\in F$, $r>0$.  We may
     assume without restriction that the diameter of $V$ is
     equal to one.  Then set $a:=\frac{\vol(V)}{\vol(B(0,2))}$, where
     $\vol$ denotes the Lebesgue measure.  Then the number of disjoint
     copies of $V$ intersecting a ball of radius $1$ is not greater
     than the number of disjoint copies of $V$ contained in a ball of
     radius $2$ which is in turn lower equal than $a^{-1}$.
	
     As a consequence, for any $x\in F$, if
     $\displaystyle\cai(x,n):=\{\s\in\S_{n}:V_{\s}\cap B_{F}(x,\La_{n}) \ne
     \emptyset\}$, then
     \begin{equation}\label{VinBall}
	 \#\cai(x,n) \leq a^{-1}.
     \end{equation}
     Clearly, by the regularity of $V$, $B_{F}(x,\La_{n}) \subseteq
     V_{\cai(x,n)}$, hence
     \begin{equation}\label{eq:muleq}
	\m(B_{F}(x,\La_{n})) \leq \m(V_{\cai(x,n)}) \leq \frac{1}{a P_{n}}.
     \end{equation}
     On the other hand, if $x\in F$, there is $\s(x)\in\S_{n}$ such 
     that $x\in\ov{V_{\s(x)}}$, therefore, for any $r>\La_{n}$, 
     $B_{F}(x,r) \supseteq C_{\s(x)}$, hence
     \begin{equation}\label{eq:mugeq}
	 \m(B_{F}(x,r)) \geq\m(C_{\s(x)}) \geq\frac{1}{P_{n}}.
     \end{equation}
     Then, for any $r>0$, if $n = n_{r}\in\bn$ is such that $\La_{n}
     < r \leq \La_{n-1}$, we get, for $x,\ y\in F$,
     \begin{align*}
	 \frac{a}{p} & \leq \frac{a}{p_{n}} =
	 \frac{1/P_{n}}{1/(aP_{n-1})} \\
	 & \leq \frac{\m(B_{F}(x,r))}{\m(B_{F}(y,\La_{n-1}))} \leq
	 \frac{\m(B_{F}(x,r))}{\m(B_{F}(y,r))} \leq
	 \frac{\m(B_{F}(x,\La_{n-1}))}{\m(B_{F}(y,r))} \\
	 & \leq \frac{1/(aP_{n-1})}{1/P_{n}} = \frac{p_{n}}{a}
	 \leq\frac{p}{a}.
     \end{align*}
     This proves Condition \ref{NewAssum}, therefore the other
     statements follow from Corollaries \ref{cor:ptw.box=ptw.mis},
     \ref{cor:new}.
 \end{proof}
 \begin{thm}\label{TgDimForLimFrac}
     Let $F$ be a translation fractal with regular open set condition, with
     the notation above, and assume $p:=\sup_{n}p_{n}<\infty$. Then
     \item{$(i)$}
     \begin{align*}
	\underline{\d}_{\m}(x) &= \liminf_{n,k\to\infty}\frac{\log
	P_{n+k}-\log P_{n}} {\log 1/\La_{n+k}-\log 1/\La_{n}},\\
	\overline{\d}_{\m}(x) &= \limsup_{n,k\to\infty}\frac{\log
	P_{n+k}-\log P_{n}} {\log 1/\La_{n+k}-\log 1/\La_{n}}.
    \end{align*}
    \item{$(ii)$}
    \begin{align*}
	\subd_{\m}(x) & = d_{H}(F) = \liminf_{n\to\infty} \frac{\log
	P_{n}}{\log 1/\La_{n}},\\
	\supd_{\m}(x) & = \limsup_{n\to\infty} \frac{\log P_{n}}{\log
	1/\La_{n}}.
    \end{align*}
    Moreover the Hausdorff measure corresponding to $d:=d_{H}(F)$ is
    non trivial if and only if $\liminf(\log P_{n}-d \log 1/\La_{n})$
    is finite.
 \end{thm} 
 \begin{proof}
     $(i)$ Let us denote by $\cl{p,\l}$ the set of limit points, for
     $n,k\to+\infty$, of $\frac{\log P_{n+k}-\log P_{n}} {\log
     1/\La_{n+k}-\log 1/\La_{n}}$, and by $\cl{f}$ the set of limit
     points, for $t,h\to+\infty$, of $\frac{f(t+h)-f(t)}{h}$, where we
     set $f(t):=-\log\m(B(x,e^{-t}))$.  Recalling Theorem \ref{eqdef},
     $(iii)$, the formulas are proved if we show that
     $\cl{p,\l}\subseteq\cl{f}$, and that for any $c\in\cl{f}$ there
     exist $c',c''\in\cl{p,\l}$ such that $c'\leq c\leq c''$. 
     Concerning the inclusion, from (\ref{eq:muleq}) and
     (\ref{eq:mugeq}) we have
     \begin{align*}
	 -\log1/a+f(\log 1/\La_{n+k}) & -f(\log 1/\La_{n}) \leq \log
	 P_{n+k}-\log P_{n} \\
	 &\leq \log1/a+f(\log 1/\La_{n+k})-f(\log 1/\La_{n}).
     \end{align*}
     As a consequence, if $n_{j}, k_{j}$ are subsequences giving 
     rise to a limit point in $\cl{p,\l}$, and we set 
     $t'_{j}=\log 1/\La_{n_{j}+k_{j}}$, $t_{j}=\log 1/\La_{n_{j}}$, 
     we obtain that $\frac{f(t'_{j})-f(t_{j})}{t'_{j}-t_{j}}$ 
     converges to the same limit, where we used that 
     $h_{j}:=\log 1/\La_{n_{j}+k_{j}}-\log 1/\La_{n_{j}}\to\infty$, since is 
     minorized by $\frac{\log 2}{N} k_{j}$.
     \\
     Now let $c\in\cl{f}$.  Then we find two sequences, $t_{k}$ and
     $t'_{k}$, such that $t_{k}\to\infty$ and
     $(t'_{k}-t_{k})\to\infty$, and
     $c=\lim_{k}\frac{f(t'_{k})-f(t_{k})}{t'_{k}-t_{k}}$.  If
     $\log1/\La_{n_{k}}$ is the best approximation from below of
     $t_{k}$ and $\log1/\La_{n'_{k}}$ is the best approximation from
     below of $t'_{k}$, we get
     \begin{align*}
	 f(t'_{k})-f(t_{k})
	 &\leq f(\log1/\La_{n'_{k}+1})-f(\log1/\La_{n_{k}})\\
	 &\leq \log P_{n'_{k}+1} - \log P_{n_{k}}-\log a\\
	 &\leq \log P_{n'_{k}} - \log P_{n_{k}+1}+2\log p-\log a,
     \end{align*}
     which shows in particular that $P_{n'_{k}} / 
     P_{n_{k}+1}\to\infty$, hence, for the bound on $p_{j}$, also 
     $n'_{k}-(n_{k}+1)\to\infty$.
     We also get
     \begin{align*}
	 t'_{k}-t_{k}
	 &\geq \log1/\La_{n'_{k}}-\log1/\La_{n_{k}+1}
     \end{align*}
     therefore
     \begin{align*}
	 \lim_{k}\frac{f(t'_{k})-f(t_{k})}{t'_{k}-t_{k}}
	 &\leq \limsup_{k}
	 \frac{\log P_{n'_{k}} - \log P_{n_{k}+1}}
	 {t'_{k}-t_{k}}\\
	 &\leq \limsup_{k}
	 \frac{\log P_{n'_{k}} - \log P_{n_{k}+1}}
	 {\log1/\La_{n'_{k}}-\log1/\La_{n_{k}+1}}.
     \end{align*}
     Possibly passing to a subsequence we obtain $c\leq 
     c''\in\cl{p,\l}$. The point $c'$ is obtained analogously.
    
     $(ii)$ Let $\{t_{k}\}$ be an increasing sequence of positive real 
    numbers s.t. $\subd_{\m}(x) = \lim_{k\to\infty} 
    \frac{f(t_{k})}{t_{k}}$, and let $\{n_{k}\}$ be an increasing sequence 
    of natural numbers s.t. $\log \frac{1}{ \La_{n_{k}} } \leq t_{k} < \log 
    \frac{1}{ \La_{ n_{k+1} } }$. Then
    \begin{align*}
	\frac{ f(t_{k}) }{ t_{k} } & \geq \frac{f( \frac{1}{
	\La_{n_{k}} } )} { \log \frac{1}{ \La_{n_{k+1}} } } \geq
	\frac{ \log P_{n_{k}} + \log a }{ \log \frac{1}{ \La_{ n_{k+1}
	} } } \\
	& \geq \frac{ \log P_{n_{k+1}} -\log p + \log a }{ \log
	\frac{1}{ \La_{ n_{k+1} } } }
    \end{align*}
    so that $\subd_{\m}(x) \geq \liminf_{n\to\infty} \frac{\log P_{n}
    }{ \log \frac{1}{ \La_{n} } }$.  Conversely, let $\{n_{k}\}$ be an
    increasing sequence of natural numbers s.t. $\lim_{k\to\infty}
    \frac{\log P_{n_{k}} }{ \log \frac{1}{ \La_{n_{k}} } } =
    \liminf_{n\to\infty} \frac{\log P_{n} }{ \log \frac{1}{ 
    \La_{n} } }$ and set $t_{k}:= \log \frac{1}{ \La_{n_{k}}}$. Then
    \begin{align*}
	\subd_{\m}(x) & = \liminf_{t\to\infty} \frac{f(t)}{t} \leq 
	\liminf_{k\to\infty} \frac{f(t_{k})}{t_{k}} \\
	& = \liminf_{k\to\infty} \frac{f( \frac{1}{\La_{n_{k}}} )}{
	\log \frac{1}{ \La_{n_{k}} } } \leq \liminf_{k\to\infty}
	\frac{ \log P_{n_{k}} }{ \log \frac{1}{\La_{ n_{k} } } }\\
	& = \liminf_{n\to\infty} \frac{\log P_{n} }{ \log \frac{1}{
	\La_{n} } }.
    \end{align*}
    The equation for $\supd_{\m}(x)$ is proved similarly.
    
    The equality $\subd_{\m}(x) = d_{H}(F)$ and the last statement
    follow from Theorem 1.8 in \cite{GuIs10}.
\end{proof}

\end{document}